\documentclass[12pt]{article}
\usepackage{amssymb,amsmath}

\newtheorem{theorem}{Theorem}[section]
\newtheorem{lemma}[theorem]{Lemma}

\newtheorem{definition}[theorem]{Definition}

\newtheorem{remark}[theorem]{Remark}
\newtheorem{example}[theorem]{Example}

\begin{document}

\title{Riesz basis for strongly continuous groups.
}

\author{ Hans Zwart\thanks{University of Twente, Faculty of Electrical Engineering, Mathematics and Computer Science, Department of Applied Mathematics, P.O. Box 217, 7500 AE Enschede, The Netherlands, \tt h.j.zwart@math.utwente.nl}
}

\maketitle

\begin{abstract}
  Given a Hilbert space and the generator of a strongly continuous group on this Hilbert space. If the eigenvalues of the generator have a uniform gap, and if the span of the corresponding eigenvectors is dense, then these eigenvectors form a Riesz basis (or unconditional basis) of the Hilbert space. Furthermore, we show that none of the conditions can be weakened. However, if the eigenvalues (counted with multiplicity) can be grouped into subsets of at most $K$ elements, and the distance between the groups is (uniformly) bounded away from zero, then the spectral projections associated to the groups form a Riesz family. This implies that if in every range of the spectral projection we construct an orthonormal basis, then the union of these bases is a Riesz basis in the Hilbert space.  
\end{abstract}

\section{Introduction and main results}

We begin by introducing some notation. By $H$ we denote the Hilbert space
with inner product $\langle \cdot,\cdot \rangle$ and norm $\|\cdot\|$,
and by $A$ we denote an unbounded operator from its domain $D(A)
\subset H$ to $H$.

If $A$ is self-adjoint and has a compact resolvent operator, then it has an orthonormal basis of eigenvectors. Unfortunately, even a slight perturbation of $A$ can destroy the self-adjointness of $A$ and so also the orthonormal basis property of the eigenvectors. However, in general the (normalized) eigenvectors, $\{\phi_n\}_{n \in {\mathbb N}}$  will still form a Riesz basis, i.e., their span is dense in $H$ and there exist (positive) constants $m$ and $M$ such that
\begin{equation}
  \label{eq:1.1}
    m \sum_{n=1}^N |\alpha_n|^2 \leq \left\|\sum_{n=1}^N \alpha_n \phi_n \right\|^2 \leq  m \sum_{n=1}^N |\alpha_n|^2
\end{equation}
for every sequence $\{\alpha_n\}_{n=1}^N$.
If $A$ processes a Riesz basis of eigenvectors, then many system theoretic properties like stability, controllability, etc.\ are easily checkable, see e.g.\ \cite{CuZw95}. 

Since this Riesz-basis property is so important there is an extensive
literature on this problem. We refer to the book of Dunford and
Schwartz \cite{DuSc71}, where this problem is treated for discrete
operators, i.e., the inverse of $A$ is compact. They apply these
results to differential operators. Riesz spectral properties for
differential operators is also the subject of Mennicken and Moller
\cite{MeMo03} and Tretter \cite{Tret00}. In these references no use is
made of the fact that for many differential operators the abstract
differential equation
\begin{equation}
  \label{eq:1.1a}
  \dot{x}(t) = A x(t), \qquad x(0)=x_0
\end{equation}
on the Hilbert space $H$ has a unique solution for every initial
condition $x_0$, i.e., $A$ is the infinitesimal generator of a
$C_0$-semigroup. In \cite{XuYu05} this property is used. The property
will also be essential in our paper. In many applications the differential
operator $A$ arises from a partial differential equation, for which it is
known that (\ref{eq:1.1a}) has a solution. Hence the assumption that $A$ generates a
$C_0$-semigroup is not strong. However, for our result the semigroup
property does not suffice, we need that $A$ generates a group, i.e.,
(\ref{eq:1.1a}) possesses a unique solution forward and backward in
time. Since we also assume that the eigenvalues lie in a strip
parallel to the imaginary axis, the group condition is not very
restrictive. For more information on groups and semigroups, we refer
the reader to \cite{CuZw95,EnNa00}.

The approach which we take to prove our result is different to the one
taken in \cite{DuSc71, MeMo03, Tret00}, and \cite{XuYu05}. We use
the fact that every generator of a group has a bounded ${\mathcal
  H}_{\infty}$-calculus on a strip. This means that to every complex
valued function $f$ bounded and analytic in a strip parallel to the
imaginary axis, there exists a bounded operator $f(A)$. For more detail
we refer to \cite{Haas06}. Note that the ${\mathcal
  H}_{\infty}$-calculus extends the functional calculus of Von~Neumann
\cite{Neum96} for self-adjoint operators.

We formulate our main results.
\begin{theorem}
  \label{T1.1}
  Let $A$ be the infinitesimal generator of the $C_0$-group $(T(t))_{t
    \in {\mathbb R}}$ on the Hilbert space $H$. We denote the
  eigenvalues of $A$ by $\lambda_n$ (counting with multiplicity), and
  the corresponding (normalized) eigenvectors by $\{\phi_n\}$. If the
  following two conditions hold,
  \begin{enumerate}
  \item The span of the eigenvectors form a dense set in $H$;
  \item The point spectrum has a {\em uniform gap}, i.e.,
    \begin{equation}
    \label{eq:1.2}
    \inf_{n \neq m} |\lambda_n - \lambda_m| > 0.
    \end{equation}
  \end{enumerate}
  Then the eigenvectors form a Riesz basis on $H$.
\end{theorem} 

\mbox{}From this result, we have some easy consequences.
\begin{remark}%
\label{R1.2}
  \begin{itemize}
  \item Since we are counting the eigenvalues with multiplicity, we see
    from (\ref{eq:1.2}) that all eigenvalues are simple.
  \item If the operator possesses a Riesz basis of eigenvectors, then
    the spectrum equals the closure of the point spectrum. Using
    (\ref{eq:1.2}), we see that the spectrum of $A$, $\sigma(A)$,
    is pure point spectrum, and $\sigma(A) =\{\lambda_n\}_{n\in
      {\mathbb N}}$.
  \item It is easy to see that for $\lambda \in \rho(A)$, the
    finite-range operator $\sum_{n=1}^N \frac{1}{\lambda-\lambda_n}
    \phi_n$ converges uniformly to $(\lambda I - A)^{-1}$. Thus under
    the conditions of Theorem \ref{T1.1}, the resolvent operator of
    $A$ is compact, i.e., $A$ is discrete.
  \end{itemize}
\end{remark}

If the eigenvalues do not satisfy (\ref{eq:1.2}), then Theorem \ref{T1.1} need not to hold. A simple counter example is given next.
\begin{example}
  \label{E1.3}  Let $H$ be a Hilbert space with orthonormal basis $\{e_n\}_{n \in {\mathbb N}}$. Define $A$ as
  \begin{equation}
    \label{eq:6}
    A \sum_{n \in {\mathbb N}} \alpha_n e_n = \sum_{k \in {\mathbb N}} (ki \alpha_{2k-1} + \alpha_{2k} ) e_{2k-1} + (k+\frac{1}{k})i \alpha_{2k} e_{2k}
  \end{equation}
  with domain
  \begin{equation}
    \label{eq:5}
    D(A) =\left\{ x= \sum_{n=1}^{\infty} \alpha_n e_n \in H \mid \sum_{n=1}^{\infty} |n\alpha_n|^2 < \infty \right\}.
  \end{equation}
  Hence the operator $A$ is block diagonal, i.e.,
  \begin{equation}
    \label{eq:7}
    A = \mathrm{diag}\, \left[ \begin{array}{cc} ki & 1 \\ 0 & (k +\frac{1}{k})i \end{array} \right].
  \end{equation}
  Using this structure, it is easy to see that $A$ generates a strongly continuous group on $H$, and that its eigenvalues are given by $\{ ki, (k+\frac{1}{k})i; k \in {\mathbb N}\}$, with the normalized eigenvectors $\phi_{2k-1} = e_{2k-1}$ and $\phi_{2k} = \frac{1}{\sqrt{k^2+1}} \left( ki e_{2k-1} + e_{2k} \right)$, $k \in {\mathbb N}$. 

  Since
  \[
     \inf_k \| i \phi_{2k-1} - \phi_{2k} \| =0
  \]
  we have that the eigenvectors do not from a Riesz basis.
\end{example}
  
  It is trivial to see that the eigenvalues in the above example can be written as the union of two sets with every subset satisfying (\ref{eq:1.2}). Furthermore, if we would group the eigenvectors as $\Phi_n = \{ \phi_{2n-1},\phi_{2n} \}$, then the (spectral) projections, $P_n$, on the span of $\Phi_n$ satisfy
  \begin{equation}
    \label{eq:8}
    \sum_{n} \|P_n x \|^2 = \|x\|^2.
  \end{equation}
This is equivalent to the fact that the projections are orthonormal. The concept of a Riesz basis is an extension of the concept of an orthonormal basis. Similarly, we can extend the concept of orthonormal projections.
\begin{definition}
  \label{D1.4}
  The family of projections $\{P_n, n \in {\mathbb N}\}$ is a Riesz family if there exists constants $m_1$ and $M_1$ such that
  \begin{equation}
    \label{eq:9}
    m_1 \|x\|^2 \leq \sum_{n} \|P_n x\|^2 \leq M_1 \|x\|^2
  \end{equation}
  for all $x \in H$.
\end{definition}

Please note that in \cite[section I.1.4]{AvIv95} the range of the projections is called a (Riesz) basis. We found it confusing with the standard definition of a Riesz basis, and therefor we use Riesz family. This concept is equivalent to Riesz basis in parenthesis, see \cite{Shka86}.
Wermer \cite{Werm54} proved the following characterization.
\begin{lemma}
\label{L1.5}
  A family of projections is a Riesz family if and only if there exists a $M_2$ such that for every subset ${\mathbb J}$ of ${\mathbb N}$ there holds
  \begin{equation}
    \label{eq:10}
    \|\sum_{n \in {\mathbb J}} P_n \| \leq M_2.
  \end{equation}
\end{lemma}

In the example we saw that we had a Riesz family of spectral projections. The following theorem states that this always hold when the eigenvalues can be decomposed in  a finite number of sets with every set satisfying (\ref{eq:1.2}).
\begin{theorem}
\label{T1.6}
  Let $A$ be the infinitesimal generator of the $C_0$-group $(T(t))_{t
    \in {\mathbb R}}$ on the Hilbert space $H$. We denote the
  eigenvalues of $A$ by $\lambda_n$ (counting with multiplicity). If the
  following two conditions hold,
  \begin{enumerate}
  \item The span of the (generalized) eigenvectors form a dense set in $H$;
  \item The eigenvalues $\{\lambda_n\}$ can be decomposed into $K$ sets, with every set having a uniform gap. 
  \end{enumerate}
  Then there are spectral projections $P_n$, $n \in {\mathbb N}$ such that
  \begin{enumerate}  
  \item $\sum_n P_n = I$, i.e., for every $x\in H$ there holds $\lim_{N\rightarrow \infty} \sum_{n=1}^N P_n x = x$;
  \item The dimension of the range of $P_n$ is at most $K$;
  \item The family of projections $\{P_n, n \in {\mathbb N}\}$ is a Riesz family.
  \end{enumerate}
\end{theorem}

So the difference with Theorem \ref{T1.1} is that we allow for non-simple eigenvalues, and that the eigenvalues may cluster. Apart from these differences, the other remarks of Remark \ref{R1.2} still hold. An operator satisfying the conditions of Theorem \ref{T1.6} has pure point spectrum and has compact resolvent.

\section{Functional calculus for groups.}
\label{sec:2}

We begin by introducing some notation. 
Since $\left(T(t)\right)_{t \in {\mathbb R}}$ is a group, there exists
a $\omega_0$ and $M_0$ such that $\|T(t)\| \leq M_0 e ^{\omega_0 |t|}$.

For $\alpha>0$ we define a strip parallel to the imaginary axis by $S_{\alpha}:=\{ s
\in {\mathbb C} \mid - \alpha < \mathrm{Re}(s) < \alpha \}$.

By ${\mathcal H}^{\infty}(S_{\alpha})$ we denote the linear space of
all functions from $S_{\alpha}$ to ${\mathbb C}$ which are analytic
and (uniformly) bounded on $S_{\alpha}$. The norm of a function in
${\mathcal H}^{\infty}(S_{\alpha})$ is given by
\begin{equation}
  \label{eq:2.1}
  \|f\|_{\infty}  = \sup_{s \in S_{\alpha}} |f(s)|.
\end{equation}
In Haase \cite{Haas04,Haas06} it is shown that the generator of the
group, $A$, has a ${\mathcal H}^{\infty}(S_{\alpha})$-calculus for
$\alpha > \omega_0$. This we explain in a little bit more detail.

Choose $\omega_1$ and $\omega$ such that $\omega_0 < \omega_1 <\alpha < \omega$. Furthermore, let $\Gamma = \gamma_1 \oplus \gamma_2$ with $\gamma_1 = -\omega_1- ir$, $\gamma_2=\omega_1 + i r$, $r \in {\mathbb R}$. For $f \in {\mathcal H}^{\infty}(S_{\alpha})$ and $x \in D(A^2)$ we define
\begin{equation}
\label{eq:2.2}
  f(A)x = \frac{1}{2\pi i} \int_{\Gamma} \frac{f(z)}{z^2- \omega^2} (z I - A)^{-1} d z \cdot (A^2 - \omega^2)x .
\end{equation}
For $\lambda \not\in \overline{S_{\alpha}}$ we have that 
\[
  \left( \frac{1}{\lambda - \cdot} \right)(A)x = (\lambda I - A)^{-1} x
\]
Furthermore, the operator defined in (\ref{eq:2.2}) extends to a bounded operator on $H$, and
\begin{equation}
  \label{eq:2.3}
  \|f(A)\| \leq c \|f\|_{\infty},
\end{equation}
with $c$ independent of $f$.

In the following lemma we show that this functional calculus behaves
like one would expect from the functional calculus of von~Neumann and
Dunford.
\begin{lemma}
  \label{L2.1}
  Let $A$ be the infinitesimal generator of a group and let $\phi_n$
  be an eigenvector for the eigenvalue $\lambda_n$. Then for every $f
  \in {\mathcal H}^{\infty}(S_{\alpha})$ there holds that
  \begin{equation}
    \label{eq:2.4}
      f(A) \phi_n = f(\lambda_n) \phi_n, \qquad n \in {\mathbb N}.
  \end{equation}
  Furthermore, if $\phi_{n,j}$ is the $j$-th generalized eigenvector for the eigenvalue $\lambda_n$, i.e., $(A-\lambda_n) \phi_{n,j} =\phi_{n,j-1}$, $j\geq 1$, with $\phi_{n,0}=\phi_n$, then
  \begin{equation}
    \label{eq:1}
     f(A) \phi_{n,j} = \sum_{m=0}^j \frac{f^{(j-m)}(\lambda_n)}{(j-m)!} \phi_{n,m}, \qquad n \in {\mathbb N},
 \end{equation}
  where $f^{(\ell)}$ denotes the $\ell$-th derivative of $f$.
\end{lemma}
  {\bf Proof}:\/ Since $\phi_n$ is an eigenvector, it is an element of
  $D(A^2)$ and so we may use equation (\ref{eq:2.2}).
  Hence
  \begin{eqnarray*}
    f(A)\phi_n &=&  \frac{1}{2\pi i} \int_{\Gamma} \frac{f(z)}{z^2- \omega^2} (z I - A)^{-1} d z \cdot (A^2 - \omega^2) \phi_n\\
    &=& \frac{1}{2\pi i} \int_{\Gamma} \frac{f(z)}{z^2- \omega^2} (z I - A)^{-1} d z\, (\lambda_n^2-\omega^2)\phi_n\\
    &=& (\lambda_n^2-\omega^2) \frac{1}{2\pi i} \int_{\Gamma} \frac{f(z)}{z^2- \omega^2} (z I - \lambda_n)^{-1}\phi_n d z.
  \end{eqnarray*}
  Since $f$ is bounded on $S_{\alpha}$, we see that the integrand
  converges quickly to zero for $|z|$ large. Hence we may apply Cauchy
  residue theorem. The only pole within the contour is $\lambda_n$,
  and so we obtain
  \[
    f(A) \phi_n = (\lambda_n^2-\omega^2) \frac{f(\lambda_n)}{\lambda_n^2- \omega^2} \phi_n = f(\lambda_n) \phi_n.
  \]
  This shows equation (\ref{eq:2.4}). We now prove the assertion for the generalized eigenvectors. 

  By induction it is easy to show that
  \begin{equation}
    \label{eq:2}
    (zI-A)^{-1} \phi_{n,j} = \sum_{m=0}^j \frac{1}{(z-\lambda_n)^{j+1-m}} \phi_{n,m}.
  \end{equation}
  Using (\ref{eq:2.2}) we have that
  \begin{align*}
    f(A) \phi_{n,j}  &=  \frac{1}{2\pi i} \int_{\Gamma} \frac{f(z)}{z^2- \omega^2} (z I - A)^{-1} d z \cdot (A^2 - \omega^2) \phi_{n,j}\\
      &= \frac{1}{2\pi i} \int_{\Gamma} \frac{f(z)}{z^2- \omega^2} (z I - A)^{-1} d z ((\lambda_n^2 - \omega^2)\phi_{n,j} + 2\lambda_n \phi_{n,j-1} + \phi_{n,j-2})\\
      &=  (\lambda_n^2-\omega^2) \frac{1}{2\pi i} \int_{\Gamma} \frac{f(z)}{z^2- \omega^2} \sum_{m=0}^j (z I - \lambda_n)^{-(j-m+1)}\phi_{n,m} d z + \\
   &\qquad 2\lambda_n \frac{1}{2\pi i} \int_{\Gamma} \frac{f(z)}{z^2- \omega^2} \sum_{m=0}^{j-1} (z I - \lambda_n)^{-(j-m)}\phi_{n,m} d z + \\
   &\qquad
  \frac{1}{2\pi i} \int_{\Gamma} \frac{f(z)}{z^2- \omega^2} \sum_{m=0}^{j-2} (z I - \lambda_n)^{-(j-m-1)}\phi_{n,m} d z \\
   &= (\lambda_n^2-\omega^2) \sum_{m=0}^{j} \frac{g^{(j-m)}(\lambda_n)}{(j-m)!} \phi_{n,m} + \\
   &\qquad 2\lambda_n \sum_{m=0}^{j-1} \frac{g^{(j-1-m)}(\lambda_n)}{(j-m-1)!} \phi_{n,m} +
   \sum_{m=0}^{j-2} \frac{g^{(j-2-m)}(\lambda_n)}{(j-m-2)!} \phi_{n,m},
  \end{align*}
  where we have introduced $g(z)=f(z)/(z^2-\omega^2)$. Using the fact that $f(z) = (z^2-\omega^2) g(z)$, we see that for $\ell \geq 2$
  \[
     f^{(\ell)}(z) = (z^2-\omega^2) g^{(\ell)}(z) + 2 z\ell g^{(\ell-1)}(z) + \ell(\ell-1) g^{(\ell-2)}(z) 
  \]
  and
  \[
     f^{(1)}(z) = (z^2-\omega^2) g^{(1)}(z) + 2 z g(z).
  \]
  This implies that
  \begin{eqnarray*}
   f(A) \phi_{n,j} &=& \sum_{m=0}^{j-2} \left[ (\lambda_n^2-\omega^2) \frac{g^{(j-m)}(\lambda_n)}{(j-m)!} +\right. \\
   && \hphantom{\sum_{m=0}^{j-2} \left[ \right]}
   \left. 2\lambda_n \frac{g^{(j-1-m)}(\lambda_n)}{(j-m-1)!} +
   \frac{g^{(j-2-m)}(\lambda_n)}{(j-m-2)!} \right] \phi_{n,m} + \\
   &&
  \left[ (\lambda_n^2 - \omega^2) g^{(1)}(\lambda_n) + 2 \lambda_n g(\lambda_n) \right] \phi_{n,j-1}+ (\lambda_n^2 - \omega^2) g(\lambda_n) \phi_{n,j}\\
  &=&  \sum_{m=0}^{j-2} \frac{1}{(j-m)!}f^{(j-m)}(\lambda_{n}) \phi_{n,m} + f^{(1)}(\lambda_n)  \phi_{n,j-1}  + f(\lambda_n) \phi_{n,j}.
  \end{eqnarray*}
  Hence we have proved the assertion. \hfill$\Box$
\medskip

\section{Interpolation sequences}

For the proof of Theorem \ref{T1.1} we need the following
interpolation result, see \cite[Theorem VII.1.1]{Garn07}.
\begin{theorem}
  \label{T2.2}
  Consider the sequence ${\mu_n}$ which satisfies $\beta_1 >
  \mathrm{Re}(\mu_n) > \beta_2 >0$ and $\inf_{n\neq m} |\mu_n - \mu_m|
  >0$. Then for every bounded sequence $\{\alpha_n\}_{n\in {\mathbb
      N}}$ of complex numbers there exists a function $g$ holomorphic
  and bounded in the right-half plane ${\mathbb C}_+:=\{ s \in
  {\mathbb C} \mid \mathrm{Re}(s) >0\}$ such that
  \begin{equation}
    \label{eq:2.5}
      g(\mu_n) = \alpha_n.
    \end{equation}
  Furthermore, there exists an $M$ independent of $g$ and $\{\alpha_n\}_{n \in {\mathbb N}}$ such that
  \begin{equation}
    \label{eq:2.6}
      \sup_{\mathrm{Re}(s) >0} |g(s)| \leq M \sup_{n \in {\mathbb N}} |\alpha_n|.
  \end{equation}
\end{theorem}

A sequence $\{\mu_n\}_{n \in {\mathbb N}}$ satisfying the conditions of the above theorem is called an {\em interpolation sequence}. It is well known that for any interpolation sequence, we can find a Blaschke product with exactly this sequence as its zero set. Using Schwartz lemma, and Lemma VII.5.3 of \cite{Garn07} the following is easy to show.
  \begin{lemma}
  \label{L3.2}
    Let $\{\mu_n\}_{n \in {\mathbb N}}$ be an interpolating sequence, and let $0< \delta := \inf_{n\neq m} |\mu_n - \mu_m|$. Furthermore, let $B(a, r)$ denote the ball in the complex plane with center $a$ and radius $r$. 

  Let $f$ be a function defined on $B(\mu_n,\delta/2)$ and which is analytic and bounded on this set. Furthermore, $f$ is zero at $\mu_n$, i.e., $f(\mu_n)=0$. Then there exists a constant $m_1$ independent of $n$ and $f$ such that
  \begin{equation}
    \label{eq:11}
      \sup_{s \in B(\mu_n,\delta/2)} \left| \frac{f(s)}{Bl(s)} \right| \leq m_1 \sup_{s \in B(\mu_n,\delta/2)} |f(s)|,
  \end{equation}
  where $Bl(s)$ is Blaschke product with zeros $\{\mu_n\}_{n \in {\mathbb N}}$.

  Similarly,  let $f$ be a function defined on $B(\mu_n,\delta/2)$ and which is analytic and bounded on this set. Furthermore, $f$ has $\ell$ zeros at $\mu_n$, i.e., $f(\mu_n)=\cdots = f^{(\ell)}(\mu_n)=0$. Then there exists a constant $m_\ell$ independent of $n$ and $f$ such that
  \begin{equation}
    \label{eq:11a}
      \sup_{s \in B(\mu_n,\delta/2)} \left| \frac{f(s)}{Bl(s)^\ell} \right| \leq m_\ell \sup_{s \in B(\mu_n,\delta/2)} |f(s)|.
  \end{equation}
  \end{lemma}
\begin{theorem}
  \label{T3.3} 
  Let $\{\mu_n\}_{n \in {\mathbb N}}$ be a sequence of numbers in the right-half plane, and let $\{\zeta_n\}_{n \in {\mathbb N}}$ be an interpolation sequence. Furthermore, let $K$ be a positive natural number. Suppose that we can find positive real numbers $r_n$ such that 
  \begin{itemize}
  \item $0< \inf_n r_n \leq \sup_n r_n < \infty$;
  \item $\{\mu_n\}_{n \in {\mathbb N}} \subset \cup_{n\in {\mathbb N}} B(\zeta_n,r_n)$;
  \item The number of $\mu_n$'s in $B(\zeta_k,r_k)$ is bounded by $K$, and 
  \item The distance between every two balls $B(\zeta_n,r_n)$ and $B(\zeta_m,r_m)$, $n \neq m$, is bounded away from zero.
  \end{itemize}
  Let $\nu$ be an element of ${\mathbb N} \cup \{0\}$. Under these assumptions we have that for every bounded sequence $\{\alpha_n\}_{n \in {\mathbb N}}$ there exists a function $g$ holomorphic and bounded in the right-half plane ${\mathbb C}_+$ such that
  \begin{equation}
    \label{eq:4}
     g(\mu_k) = \alpha_n \quad \mbox{if}\quad  \mu_k \in B(\zeta_n,r_n)
  \end{equation}
  and 
  \begin{equation}
    \label{eq:4a}
     g^{(q)}(\mu_k) = 0 \quad \mbox{for}\quad  1\leq q \leq \nu,\mbox{ and } \mu_k \in B(\zeta_n,r_n)
  \end{equation}
  Furthermore, there exists an $M$ independent of $g$ and $\{\alpha_n\}_{n \in {\mathbb N}}$ such that
  \begin{equation}
    \label{eq:16}
    \sup_{\mathrm{Re}(s) >0} |g(s)| \leq M \sup_{n \in {\mathbb N}} |\alpha_n|.
  \end{equation}
  Hence the function $g$ interpolates the points $\mu_k$, but points close to each other are given the same value.
\end{theorem}
  {\bf Proof}:\/ We present the proof for the case that $\kappa =2$ and $\nu=1$. For higher values of $\kappa$ and $\nu$ the proof goes similarly. We denote the points $\mu_n \in B(\zeta_n,r_n)$ by $\xi_n$ and $\gamma_n$, and we assume that $\gamma_n \neq \xi_n$.  By the assumptions, we know that these are interpolation sequences.

Let $Bl_1(s)$ and $BL_2(s)$ be the Blaschke products with zeros $\{\xi_n\}_{n \in {\mathbb N}}$ and $\{\gamma_n\}_{n \in {\mathbb N}}$, respectively.

We write $g$ in the following form
\begin{equation}
  \label{eq:12}
  g(s) = g_1(s) + Bl_1(s)g_2(s) + Bl_1(s)^2 g_3(s) + Bl_1(s)^2Bl_2(s) g_4(s), 
\end{equation}
and we construct bounded analytic functions $g_j$ such that (\ref{eq:4})--(\ref{eq:16}) are satisfied. For our situation, the interpolation conditions become $g(\xi_n)=g(\gamma_n)=\alpha_n$ and $g^{(1)}(\xi_n)=g^{(1)}(\gamma_n)=0$. Using the form (\ref{eq:12}), these conditions are equivalent to
\begin{align}
  \label{eq:3}
   g_1(\xi_n) = & \, \alpha_n\\ 
   \label{eq:13}
   g_2(\xi_n) = & -\frac{g_1^{(1)}(\xi_n)}{Bl_1^{(1)}(\xi_n)}\\
   \label{eq:14}
   g_3(\gamma_n) =& \,\frac{\alpha_n - g_1(\gamma_n)- Bl_1(\gamma_n) g_2(\gamma_n)}{Bl_1(\gamma_n)^2} \\
   \label{eq:15}
   g_4(\gamma_n) = &  \left[ g_1^{(1)}(\gamma_n) + Bl_1(\gamma_n)g_1^{(1)}(\gamma_n) + Bl_1^{(1)}(\gamma_n) g_2(\gamma_n) +\right.\\
\nonumber
    &\left. 2 Bl_1(\gamma_n) Bl_1^{(1)}(\gamma_n)g_3(\gamma_n) + Bl_1(\gamma_n)^2g_3^{(1)}(\gamma_n)\right]\left[Bl_1(\gamma_n)^2Bl_2^{(1)}(\gamma_n)\right]^{-1}.
\end{align}
Since $\{\xi_n\}$ is an interpolation sequences, and since $\{\alpha_n\}$ is a bounded sequence, we have by Theorem \ref{T2.2} the existence of a bounded $g_1$ for which (\ref{eq:2.6}) holds.

\mbox{} Define the function $f(s)=-g_1(s)+\alpha_n$. This is clearly bounded and analytic on $B(\zeta_n,r_n)$ and zero at $s=\xi_n$. By Lemma \ref{L3.2} the function $f(s)/Bl_1(s)$ is bounded (uniformly in $n$) in this ball. l'Hopital gives that the value at $s=\xi_n$ equals the right-hand side of (\ref{eq:13}), and so the right hand-side is a bounded sequence. Furthermore, since $\{\xi_n\}_{n \in {\mathbb N}}$ is an interpolation sequence we can apply Theorem \ref{T2.2} and find a bounded analytic function $g_2$ satisfying (\ref{eq:13}).

Consider the function
\begin{equation}
  \label{eq:17}
  q_1(s) = \frac{\alpha_n - g_1(s)- Bl_1(s) g_2(s)}{Bl_1(s)^2} 
\end{equation}
By (\ref{eq:3}) and (\ref{eq:13}), we have that the denominator has two zero's in $\xi_n$. Furthermore, it is analytic on $B(\zeta_n,r_n)$ and bounded independently of $n$. Lemma \ref{L3.2} implies that $q_1(s)$ is uniformly bounded. In particular, the sequence $q_1(\gamma_n)$ is uniformly bounded. By Theorem \ref{T2.2} we construct a bounded $g_3$ satisfying (\ref{eq:14}).

It remains to show that the right hand-side of (\ref{eq:15}) is a uniformly bounded sequence. The sequence $-\frac{g_3^{(1)}(\gamma_n)}{Bl_2^{(1)}(\gamma_n)}$ is uniformly bounded for the same reason why the sequence in (\ref{eq:13}) was bounded. So we may disregard that term in (\ref{eq:15}). Using the value of $g_3(\gamma_n)$ as found in (\ref{eq:14}), we have that the denominator of (\ref{eq:15}) becomes
\begin{align}
  \label{eq:18}
    g_1^{(1)}(\gamma_n) + Bl_1(\gamma_n)g_1^{(1)}(\gamma_n) +& Bl_1^{(1)}(\gamma_n) g_2(\gamma_n) + \\
\nonumber
   &2 Bl^{(1)}_1(\gamma_n) \frac{\alpha_n - g_1(\gamma_n)- Bl_1(\gamma_n) g_2(\gamma_n)}{Bl_1(\gamma_n)}.
\end{align}
Based on this and using equation (\ref{eq:17}) we define in $B(\zeta_n,r_n)$ the bounded analytic function
\begin{equation}
  \label{eq:20}
  f(s)= g_1^{(1)}(s) + Bl_1(s)g_1^{(1)}(s) + Bl_1^{(1)}(s) g_2(s) + 2 Bl_1(s)Bl^{(1)}_1(s) q_1(s).
\end{equation}
Using (\ref{eq:13}), we have that $f(\xi_n)=0$, and using that 
\[
   q_1(\xi_n) = \frac{g_1^{(2)}(\xi_n) - Bl_1^{(2)}(\xi_n) g_2(\xi_n) - 2 Bl_1^{(1)}(\xi_n) g_2^{(1)}(\xi_n)}{2 \left( Bl_1^{(1)}(\xi_n)\right)^2}
\]
we find that $f^{(1)}(\xi_n)$ is zero as well. By Lemma \ref{L3.2}, we conclude that 
  \begin{align*}
    \left[ g_1^{(1)}(\gamma_n) + Bl_1(\gamma_n) \right.  & g_1^{(1)}(\gamma_n) + Bl_1^{(1)}(\gamma_n) g_2(\gamma_n) +\\
     & \left. 2 Bl_1(\gamma_n) Bl_1^{(1)}(\gamma_n)g_3(\gamma_n) + Bl_1(\gamma_n)^2g_3^{(1)}(\gamma_n)\right]\left[Bl_1(\gamma_n)^2\right]^{-1}
  \end{align*}
  is uniformly bounded. Since $Bl_2(s)$ is a Blaschke product with zeros $\{\gamma_n\}$ we have that $Bl_2^{(1)}(\gamma_n)$ is bounded away from zero. Hence the right hand-side of (\ref{eq:15}) is a bounded sequence, and so by Theorem \ref{T2.2} we can find the interpolating $g_4$. This concludes the construction.
 \hfill$\Box$
\medskip

Now we have all the ingredients for the proof of Theorem \ref{T1.1} and \ref{T1.6}.

\section{Proof of Theorem \ref{T1.1} and \ref{T1.6} }

As may be clear from the formulation of the Theorems \ref{T1.1} and \ref{T1.6}, Theorem \ref{T1.1} is a special case of Theorem \ref{T1.1}. However, since the Riesz basis property is for applications more important than Riesz family, we we decided to formulate them separately. The proof of Theorem \ref{T1.1} is more simple than that of the general theorem, but the underlying ideas are the same. 
\medskip

\noindent{\bf Proof of Theorem \ref{T1.1}}: 
Let $\alpha$ be the positive number defined at the beginning of Section \ref{sec:2}. We define the complex numbers $\mu_n$ as
\begin{equation}
  \label{eq:2.7}
  \mu_n = \lambda_n + \alpha \qquad n \in {\mathbb N}.
\end{equation}
By the conditions on $\alpha$ and $\lambda_n$, we see that $\{\mu_n\}_{n \in {\mathbb N}}$ satisfies the conditions of Theorem \ref{T2.2}.

Let ${\mathbb J}$ be a subset of ${\mathbb N}$. Since the eigenvalues
satisfy (\ref{eq:1.2}), we conclude by Theorem \ref{T2.2} there exists
a function $g_{{\mathbb J}}$ bounded and analytic in ${\mathbb C}_+$
such that
\begin{equation}
  \label{eq:2.8}
  g_{\mathbb J}(\mu_n)= 
  \begin{cases} 1, & \quad \mbox{if } n \in {\mathbb J}\\
                0, & \quad \mbox{if } n \not\in {\mathbb J}. 
  \end{cases} 
\end{equation}
Furthermore, see (\ref{eq:2.6})
\begin{equation}
  \label{eq:2.10}
  \sup_{s \in {\mathbb C}_+} |g_{{\mathbb J}}(s)| \leq M.
\end{equation}

Given this $g_{\mathbb J}$ we define $f_{\mathbb J}$ as
\begin{equation}
  \label{eq:2.11}
  f_{\mathbb J}(s) = g_{\mathbb J}(s+ \alpha), \qquad s \in S_{\alpha}.
\end{equation}
Then using the properties of $g_{\mathbb J}$ we have that $f_{\mathbb J} \in {\mathcal H}_{\infty}(S_{\alpha})$, 
\begin{equation}
  \label{eq:2.12}
  f_{\mathbb J}(\lambda_n)= 
  \begin{cases}  1, & \quad \mbox{if } n \in {\mathbb J}\\
                 0, & \quad \mbox{if } n \not\in {\mathbb J},
  \end{cases} 
\end{equation}
and there exists a $M>0$ independent of $f_{\mathbb J}$ such that 
\begin{equation}
  \label{eq:2.13}
  \|f_{{\mathbb J}}\|_{\infty} \leq M.
\end{equation}

Next we identify the operator $f_{\mathbb J}(A)$. Combining
(\ref{eq:2.4}) with (\ref{eq:2.12}) gives
\[
  f_{\mathbb J}(A)\phi_n = 
  \begin{cases}  \phi_n, & \quad \mbox{if } n \in {\mathbb J}\\
                 0, & \quad \mbox{if } n \not\in {\mathbb J}.
  \end{cases}
\]
Since $f_{\mathbb J}(A)$ is a linear operator, we obtain that
\begin{equation}
  \label{eq:2.14}
  f_{\mathbb J}(A) \left(\sum_{n=1}^N \alpha_n \phi_n\right) = \sum_{n\in {\mathbb J} \cap \{1,\cdots,N\}} \alpha_n \phi_n.
\end{equation}
By assumption the span of $\{\phi_n\}_{n \in {\mathbb N}}$ is dense in
$H$. Furthermore, $f_{\mathbb J}(A)$ is a bounded operator. So we
conclude that $f_{\mathbb J}$ is the spectral projection associated to
the spectral set $\{ \lambda_n \mid n \in {\mathbb J}\}$.

Combining (\ref{eq:2.3}) with (\ref{eq:2.13}) we have that these
spectral projections are uniformly bounded. Since the eigenvalues are
simple, this implies that the (normalized) eigenvectors form a Riesz
basis, see Lemma \ref{L1.5}. \hfill$\Box$
\medskip

\noindent{\bf Proof of Theorem \ref{T1.6}}:
As in the previous proof we can shift the eigenvalues by $\alpha$ such that they all lies in the right half plane, and they are bounded away from the imaginary axis. We denote these shifted eigenvalues by $\mu$. 
Since the $\lambda_n$'s can be decomposed into $K$ interpolation sequences, the same holds for $\mu_n$. Hence we can group the $\mu_n$'s as in Theorem \ref{T3.3}. Note that we are counting the eigenvalues $\lambda_n$'s, and thus $\mu_n$, with their multiplicity, and so in every ball there can at most be $K$ different values, and the multiplicity of every value is also bounded by $K$. Let us renumber the $\mu_n$'s such that the values in the $n$'th ball are given by $\mu_{n,k}$, $k=1,\cdots,n_k$.  The eigenvectors corresponding to $\lambda_{n,k}=\mu_{n,k}-\alpha$ are denoted by $\phi_{n,k,0}$, and the generalized eigenvectors by $\phi_{n,k,j}$, $j =1,\cdots, j_{n,k}$. By construction we have that 
\begin{equation}
  \label{eq:21}
  \sum_{k=1}^{n_k} \left[ j_{n,k}+1\right] \leq K.
\end{equation}

Let ${\mathbb J}$ be a subset of ${\mathbb N}$. By the above, we conclude from Theorem \ref{T3.3} that there exists a function $g_{\mathbb J}$ bounded and analytic in ${\mathbb C}_+$ such that
\begin{equation}
  \label{eq:2.15}
  g_{\mathbb J}(\mu_{n,k})= 
  \begin{cases} 1, & \quad \mbox{if } n \in {\mathbb J}\\
                0, & \quad \mbox{if } n \not\in {\mathbb J}. 
  \end{cases} 
\end{equation}
and
\begin{equation}
  \label{eq:2.16}
   g_{\mathbb J}^{(j)}(\mu_{n,k})=0, \qquad j=1,\cdots, K.
\end{equation}
Furthermore, see (\ref{eq:16})
\begin{equation}
  \label{eq:2.17}
  \sup_{s \in {\mathbb C}_+} |g_{{\mathbb J}}(s)| \leq M.
\end{equation}

Given this $g_{\mathbb J}$ we define $f_{\mathbb J}$ as
\begin{equation}
  \label{eq:2.18}
  f_{\mathbb J}(s) = g_{\mathbb J}(s+ \alpha), \qquad s \in S_{\alpha}.
\end{equation}
Then using the properties of $g_{\mathbb J}$ we have that $f_{\mathbb J} \in {\mathcal H}_{\infty}(S_{\alpha})$, 
\begin{equation}
  \label{eq:2.19}
  f_{\mathbb J}(\lambda_{n,k})= 
  \begin{cases}  1, & \quad \mbox{if } n \in {\mathbb J}\\
                 0, & \quad \mbox{if } n \not\in {\mathbb J},
  \end{cases} 
\end{equation}
\begin{equation}
  \label{eq:2.20}
   f_{\mathbb J}^{(j)}(\lambda_{n,k})=0, \qquad j=1,\cdots, K.
\end{equation}
and there exists a $M>0$ independent of ${\mathbb J}$ such that 
\begin{equation}
  \label{eq:2.21}
  \|f_{{\mathbb J}}\|_{\infty} \leq M.
\end{equation}

Next we identify the operator $f_{\mathbb J}(A)$. Combining
(\ref{eq:1}) with (\ref{eq:2.19}) and (\ref{eq:2.20}), gives
\[
  f_{\mathbb J}(A)\phi_{n,k,j} = 
  \begin{cases}  \phi_{n,k,j}, & \quad \mbox{if } n \in {\mathbb J}\\
                 0, & \quad \mbox{if } n \not\in {\mathbb J}.
  \end{cases}
\]
Since $f_{\mathbb J}(A)$ is a linear operator, we obtain that
\begin{equation}
  \label{eq:2.22}
  f_{\mathbb J}(A) \left(\sum_{n=1}^N \sum_{k=1}^{n_k} \sum_{j=0}^{j_{n,k}} \alpha_{n,k,j} \phi_{n,k,j}\right) = \sum_{n\in {\mathbb J} \cap \{1,\cdots,N\}}\sum_{k=1}^{n_k} \sum_{j=0}^{j_{n,k}}  \alpha_{n,k,j} \phi_{n,k,j}.
\end{equation}
By assumption the span of $\{\phi_{n,k,j}\}_{n \in {\mathbb N}, k=1,\cdots,n_k,,j=0,\cdots,j_{n,k}}$ is dense in
$H$. Furthermore, $f_{\mathbb J}(A)$ is a bounded operator. So we
conclude that $f_{\mathbb J}$ is the spectral projection associated to
the spectral set $\{ \lambda_{n,k} \mid n \in {\mathbb J}\}$.

Combining (\ref{eq:2.21}) with (\ref{eq:2.22}) we have that these
spectral projections are uniformly bounded. By Lemma \ref{L1.5} we conclude that the spectral projections $\{P_n\}_{n \in {\mathbb N}}$, where $P_n$ is the spectral projection associated to the eigenvalues in the $n$-th ball, are a spectral family. From (\ref{eq:21}) we see that the dimension of the range of $P_n$ is bounded by $K$. \hfill$\Box$
\medskip

\section{Closing remarks}
\label{sec:3}

A natural question is the following: {\em If $A$ is the infinitesimal generator of a group and $A$ has only point spectrum with multiplicity one, is the span over all eigenvectors dense in $H$}? 

In general the answer to this question is negative. On page 665 of Hille and Phillips \cite{HiPh57} one may find an example of a generator of a group without any spectrum. 

However, there are some interesting cases for which the answer is positive. If $A$ generates a bounded group, i.e., $\sup_{t \in {\mathbb R}}\|T(t)\| < \infty$, then $A$ is similar to a skew-adjoint operator, and so there is a complete spectral measure, see \cite{Cast83,Zwar01}. Another interesting situation is the following. Since $A$ generates a group, it can be written as $A=A_0+Q$, where $A_0$ generates a bounded group, and $Q$ is a bounded linear operator, see \cite{Haas04}. If $A_0$ has only point spectrum which satisfies (\ref{eq:1.2}), then by Theorem XIX.5.7 of \cite{DuSc71} we know that condition 1.\ holds for $A$.

When calculating the eigenvalues of a differential operator, one normally finds that these eigenvalues are the zeros of an entire function. If this function has its zeros in a strip parallel to the imaginary axis, and on the boundary of this strip the function is bounded and bounded away from zero, then its zeros can be decomposed into finitely many interpolation sequences, see Proposition II.1.28 of \cite{AvIv95}. They name this class of entire function {\em sine type functions}, but in Levin \cite{Levi96} this name is restricted to a smaller class of functions.

\subsection*{Acknowledgment}

The author wants to thank Markus Haase and Jonathan Partington for their help .

\end{document}